%% file: EmpiricalMagnitudeXXX.tex
\documentclass[a4paper,10pt]{amsart}
\usepackage{hyperref}
\usepackage{afterpage}

\usepackage{tikz}
\usepackage{pgfplots}

\usepackage[T1]{fontenc}
\usepackage[sc]{mathpazo}

\usepackage{fancyvrb}
\DefineVerbatimEnvironment{code}{Verbatim}{fontsize=\small}

\hypersetup{
  colorlinks = false,
  urlcolor = black,
  urlbordercolor = yellow,
  linkcolor = {1,1,1},
  linkbordercolor=yellow,
  pdfpagemode = UseNone
}

\input{SWi}

\title[Calculations for the magnitude of metric spaces]{Heuristic and computer calculations\\for the magnitude of metric spaces}
\author{Simon Willerton}
\email{s.willerton@sheffield.ac.uk}
\address{School of Mathematics and Statistics, Hicks Building, University of Sheffield, Sheffield, S3 7RH, UK}

\begin{document}

\begin{abstract}
The notion of the magnitude of a compact metric space was considered in a previous paper with Tom Leinster, where the magnitude was calculated for line segments, circles and Cantor sets.  
In this paper more evidence is presented for a conjectured relationship with a geometric measure theoretic valuation.  Firstly, a heuristic is given for deriving this valuation by considering `large' subspaces of Euclidean space and, secondly, numerical approximations to the magnitude are calculated for squares, disks, cubes, annuli, tori and Sierpinski gaskets.  The valuation is seen to be very close to the magnitude for the convex spaces considered and is seen to be `asymptotically' close for some other spaces. 
\end{abstract}
 
\maketitle

\section*{Introduction}
This paper is concerned with an empirical look at the magnitude of subspaces of Euclidean space.  Magnitude was introduced under the name `cardinality' by Leinster in~\cite{Leinster:MetricSpacesNCatCafe} as a partially defined invariant of finite metric spaces; his motivation came from enriched category theory, but the invariant had already been considered in the biological diversity literature \cite{SolowPolasky:MeasuringBiologicalDiversity}.  In~\cite{LeinsterWillerton:AsymptoticMagnitude} we started to extend the definition to compact, non-finite metric spaces by approximating these with finite metric spaces; we considered, in particular, the magnitude of certain subsets of Euclidean space, such as line segments and circles and we showed that in these cases, when the metric space is scaled up in size, the magnitude asymptotically satisfies the inclusion-exclusion principle.  This leads to some natural conjectures and the current paper provides numerical and heuristic evidence for these conjectures.  The current paper, however, can be read independently of~\cite{LeinsterWillerton:AsymptoticMagnitude}.

  As hinted by the title, there are two main ingredients to this paper.  The first ingredient is a heuristic argument to calculate the contribution to the magnitude from the `bulk' of points in the closure of a `large' open subset $X$ of Euclidean space $\R^n$; this contribution is shown to be roughly $\vol(X)/n!\,\omega_n$ where $\omega_n$ is the volume of the unit $n$-ball.  From here one is led to consider, for naturality reasons, the so-called penguin valuation, defined on certain compact subsets of $\R^n$ by
     \[P(X):=\sum_{i=0}^n \frac{\mu_i(X)}{i!\,\omega_i}\]
  where $\mu_i(X)$ is the $i$th intrinsic volume of $X$ (see \cite{KlainRota:Book,Federer:CurvatureMeasures}).   The top intrinsic volume  $\mu_{n}(X)$ is just the usual $n$-dimensional volume of $X$, whereas $\mu_{n-1}(X)$ is half the `surface area', meaning half the $n-1$-dimensional volume of the boundary of $X$, and $\mu_{0}(X)$ is the Euler characteristic of $X$.  These intrinsic volumes have the scaling property that for $t>0$ if $tX$ denotes $X$ with all the distances scaled up by $t$ then $\mu_{i}(tX)=t^{i}\mu_{i}(X)$.  Furthermore they are normalized on unit cubes, so that $\mu_{i}([0,1]^{i})=1$.  I am not sure the exact class of subsets of Euclidean space on which the intrinsic volumes are defined, but it includes polyconvex sets, meaning finite unions of convex sets, \cite{KlainRota:Book} and smooth submanifolds \cite{Federer:CurvatureMeasures}.
   
   The second ingredient consists of computer calculations of finite spaces approximating squares, cubes, discs, annuli, tori and Sierpinski gaskets.  The computer algebra \texttt{maple} was used to calculate the magnitude of finite spaces with around $25,000$ points.

 The two ingredients are combined when the calculated magnitudes for the first five types of spaces mentioned above are plotted together with the penguin valuation.   In the convex examples --- squares, cubes and discs --- there is a surprisingly good fit, leading to the conjecture that for $K$ a convex subset of Euclidean space the magnitude and penguin valuation are equal: $|K|=P(K)$.  In the case that $K$ is a straight line segment this is proved rigourously in~\cite{LeinsterWillerton:AsymptoticMagnitude}.  In the non-convex examples
 there is a good fit for larger spaces, leading to the conjecture that in those cases the magnitude is asymptotically the same as the penguin functional in the sense that
  \[|tX|-P(tX)\to 0 \qquad\text{as } t\to \infty,\]
although it is not clear what spaces $X$ this would be expected to hold for in general, other than for the closure of open sets.  For instance, for the `bent line', in which two straight line segments are stuck together at an angle, numerical calculations of the magnitude by David Speyer seem to show a small, constant discrepancy with the penguin valuation;  also, asymptotic calculations seem to show that the lower order terms are not right in the case of the $3$-sphere, although I do not have much confidence in these calculations at the moment.

In the fractal case it is not possible to make sense of the penguin functional. However, in  \cite{LeinsterWillerton:AsymptoticMagnitude} we were able to calculate precisely the magnitude of the Cantor set and observed that the asymptotic growth rate was the Hausdorff dimension of the Cantor set.  In Section~\ref{Section:Sierpinski} we see that the empirical data for the Sierpinski gasket is consistent with its magnitude of the growing like its side length to the power of its Hausdorff dimension.

The rest of this introduction consists of a reminder of the definition of magnitude and an informal analogy which might provide some intuition, followed by a more detailed synopsis of the paper.

\subsection*{Definition of magnitude and a useful analogy}  We first define the magnitude on finite metric spaces, then after some comments on the definition we will say how we wish to extend it to non-finite metric spaces.

If $X$ is a finite metric space then a weighting on $X$ is an assignment of a real number $w_{x}\in\R$ to each point $x\in X$ such that for each point $x'\in X$ the \emph{weight equation} for $x'$ is satisfied:
 \[\sum_{x\in X}e^{-d(x,x')}w_{x}=1.\]
If a weighting exists then the \emph{magnitude} $\left| X\right|$ is defined to be the sum of the weights: $\left| X\right|:=\sum_{x\in X}w_{x}$.

Two standard observations to make here are, firstly, that the weights are not necessarily positive, as will be seen in examples below, and, secondly, that if there is more than one weighting then the magnitude is independent of the weighting.

One informal analogy that might be useful to keep in mind is the following.
Imagine that the points are certain organisms that wish to maintain a certain body temperature which we normalize to be a temperature of $1$ unit.  To achieve this body temperature each organism can generate heat or cold and experiences the heat or cold from other organisms in a way that falls off exponentially with the distance.  If the amount of heat that organism $x$ is generating is given by $w_x$ then the fact that an organism $x'$ maintains unit body temperature is expressed in the weight equation for $x'$.

Weightings of some finite subsets of $\R^2$ (with the subspace metric) are given in Figures~\ref{Fig:WeightsOfSquaresl2} and~\ref{Fig:SquaresVaryingWidths}:  the area of each disk is proportional to the weight of the point at its centre and the colour of the disk represents the sign of the weight --- red is positive and blue is negative.  Note how the weights are greater around the edge of the metric space.  In line with the intuition described above, one is led to think about emperor penguins which huddle together in large groups on the Antarctic ice during the winter in order to conserve heat (see~\cite{Attenborough:LifeInTheFreezerPenguins}); the penguins take turns to be on the outside of the group where it is coldest and they need to use up more energy.  Despite being cute and whimsical, this imagery does seem to give helpful intuition.

We are interested here in extending the notion of magnitude to compact subsets of Euclidean space.  If $X$ is such a compact subset then we approximate $X$ by a sequence $\{X_k\}_{k=0}^{\infty}$ of finite subsets of $X$ converging (in the Hausdorff topology) to $X$.  If $\lim_{k\to\infty}|X_k|$ exists then we define $|X|$ the magnitude of $X$ to be this limit.  In general we do not know that this is a well defined process --- maybe the limit does not exist, or the limit depends on the choice of approximating sequence, or the approximating subsets do not themselves have a magnitude --- but it gives sensible answers in the examples that we can calculate.  In this paper typically we start with a compact subset $X$ of $\R^n$ and intersect it with some lattice $\LLL$ to give a finite set $\ddot X:=X\cap \LLL$.  (This is pronounced ``$X$ dots'' and is supposed to indicate the fact that it is a discrete set.)  We assume that for a sufficiently fine lattice this gives rise to a good approximation for the magnitude we are interested in: $|\ddot X| \simeq |X|$.

\subsection*{What is in this paper}
In Section~\ref{Section:Squares} we have a look at some features of weighting numerically calculated on grid approximations to squares.  In Section~2 part of these features are fed into a back-of-the-envelope calculation for the contribution to the magnitude from the bulk of point in a large subset of $\R^n$; this calculation then leads to the definition of the penguin valuation.  Section~3 contains the main results of the paper with the penguin valuation being graphed together with numerically calculated approximations to the magnitude for some subsets of Euclidean space.  In the convex examples --- squares, discs and cubes --- the penguin valuation and the magnitude are seen to be very close; in the other examples --- annuli and tori --- they are seen to be very close when the spaces are large.  In Section~4 the Sierpinski gasket is considered and it is seen that the emprirical data is compatible with the growth rate of the magnitude being the Hausdorff dimension.  The appendix contains some comments on the computer calculations.

\subsection*{Acknowledgements}
I would like to thank Tom Leinster for many stimulating conversations and the information about emperor penguins, David Speyer for useful comments on an earlier draft, Sam Marsh for encouraging me to calculate some examples and also Emma McCabe and Richard Thomas for helping me develop the heat analogy.

\section{Numerically approximating the magnitude of a square}
\label{Section:Squares}
In this section we have our first look at the results of calculations of the weights of points on a square grid and make some observations about these which will feed into the Bulk Approximation in the next section.

Having looked at line segments, circles and Cantor sets in~\cite{LeinsterWillerton:AsymptoticMagnitude}, the next most obvious class of spaces to look at was squares.  Unlike the simpler examples, we are not yet able to find an exact formula for weighting in this case, but we can look at the results of computer calculations.  For $t>0$, let $tQ$ be a $t\times t$ square, thought of as a subset of $\R^2$, this will be approximated by $t\ddot Q$ which is a square grid of $m\times m$ points with the distance between adjacent points being $(m-1)/t$.  Remember that the magnitude $|t\ddot Q|$ is supposed to give an approximation to the magnitude $|tQ|$ of the continuous square---however, it isn't yet proved  that the continuous square has a well-defined magnitude, although the data in this paper provides evidence for the existence of such a magnitude.

For specific values of $t$ and $m$, provided that $m$ is not too large, we can use a computer algebra package such as \texttt{maple} to numerically calculate the weights on the finite metric space $t\ddot Q$, and hence calculate the magnitude $|t\ddot Q|$, which is just the sum of the weights; details of the computer implementation are given in the appendix.

\begin{figure}[thb]
\begin{center}
\includegraphics[width=0.95\textwidth]{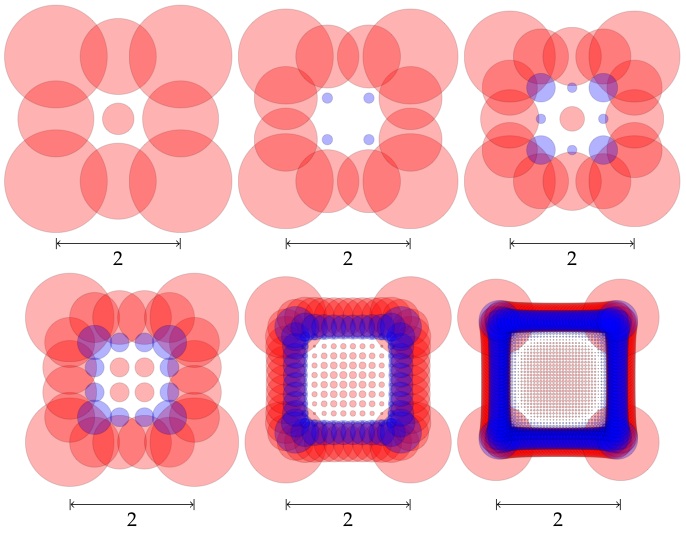}
\end{center}
\caption{Weightings on square grids of side length $2$ with respectively $3$, $4$, $5$, $6$, $14$ and $40$ points per side.
}
\label{Fig:WeightsOfSquaresl2}
\end{figure}

\begin{figure}[p]
\begin{center}
\includegraphics[width=0.85\textwidth]{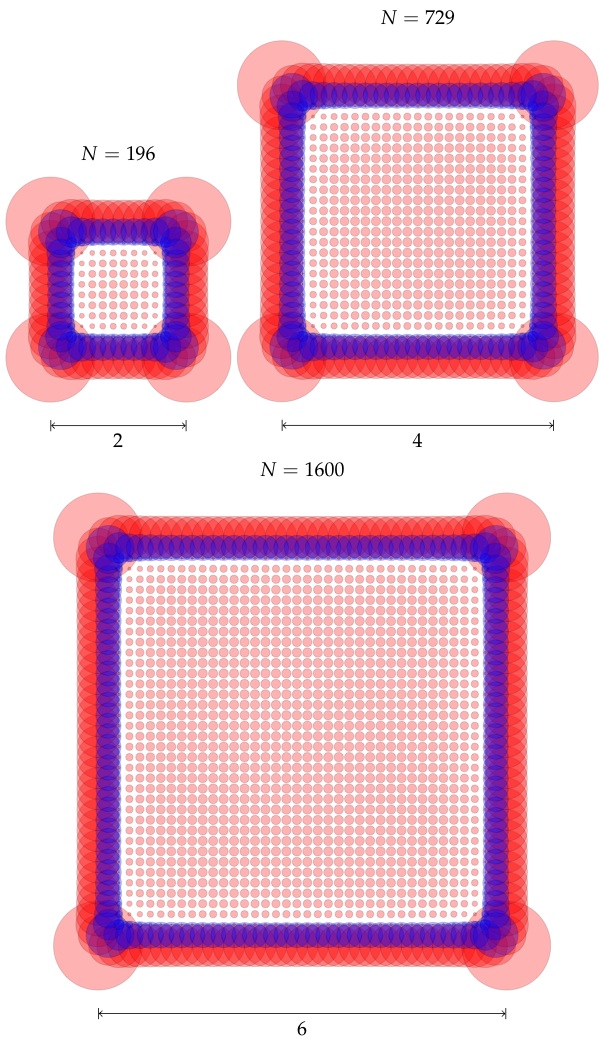}
\end{center}
\caption{Weightings on the square grid of  side length $2$, $4$ and $6$.  The grid spacing is the same, approximately $0.15$, in each case.
  The area of a disc is four
times the weight of the point at its centre; red denotes a positive weight and blue a negative weight.
}
\label{Fig:SquaresVaryingWidths}
\end{figure}

The results of some computer calculations of the weights of $t\ddot Q$ for some values of side length $t$ and $m$ points per side are represented graphically in Figures~\ref{Fig:WeightsOfSquaresl2} and~\ref{Fig:SquaresVaryingWidths}.  In these pictures the weight at each point is represented by a disc centred at that point with the area of the disc equal to four times the weight at the point; the colour of the disc represents the sign of the weight, with red meaning a positive (or `hot') weight and blue meaning a negative (or `cold') weight.
\afterpage{\clearpage}

\begin{figure}[ht]
\begin{center}
\begin{tikzpicture}
\begin{axis}[tick label style={font=\footnotesize},
axis x line=center, axis y line = left,
xmin=0,ymin=-2, ymax=1, xmax=2,
xlabel={$d$}, x label style={at={(0.9,0.52)}},
x axis line style={style = -},y axis line style={style = -}]
\addplot[only marks,mark=asterisk] file {maple/SqProfile10N171.table};
\end{axis}
\end{tikzpicture}
\end{center}
\caption{The bulk-normalized weights of the middle row of points at a distance $d$ from the edge on a $10\times 10$ square with $171\times 171$ grid of points.}
\label{Figure:ProfileNearEdge}
\end{figure}
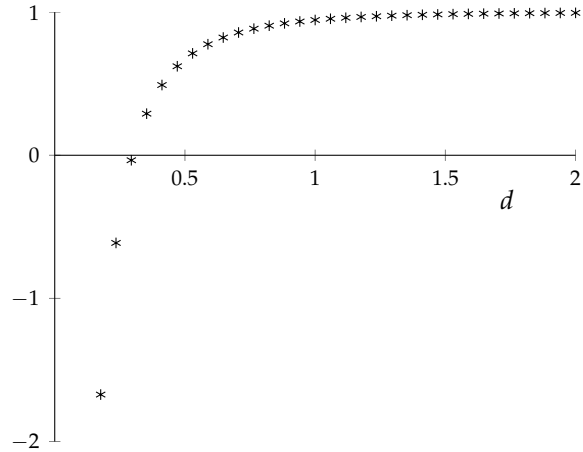

This graphical representation of the weights is not perfect: for instance where the discs overlap near the edges it is not entirely clear what is happening (see Figure~\ref{Figure:ProfileNearEdge} for a better picture of the behaviour near the edge).  However, these pictures make the following features apparent.
\begin{enumerate}
\item In the `bulk' of the square the weights are reasonably constant.
\item On the boundary the weights are bigger, but again are reasonably constant.
\item Near to the boundary the weights are negative.
\end{enumerate}
In the next section we give a heuristic argument for the first feature and this argument implies that the constant should be $1/2\pi$ times the area of the little squares in the grid.  A vague heuristic for the second and third features was given in the introduction --- you should think of penguins huddling together to keep warm, those on the outside of the huddle need to generate more heat to keep warm whilst those next to them are overheated and so need to cool themselves down.

We finish this section by looking a little closer at the behaviour near the edge.  We consider a square with $(2s+1)\times(2s+1)$ points then we look at the weights near the edge along the middle row --- that is the $s$th row.  If a point $x$ has weight $w_x$ then, based on the calculations of the next section, we define the \emph{bulk-normalized weight} of a point $x$ in the square to be $2\pi w_x/\vol(V)$ where $\vol(V)$ is the volume of a little square in the grid --- that is the volume of a Voronoi cell in the language of the next section.  According to my rather limited calculations, if the side length $t$ of the square is sufficiently large, say at least $10$, then the bulk-normalized weight of a point on the middle row near the edge depends only on the distance from the edge and not on the side length or the spacing between points, with the singular exception of the weight of the point \emph{nearest} to the edge which is usually more negative than might be expected.  A sample profile of bulk-normalized weights near the edge is given in Figure~\ref{Figure:ProfileNearEdge}.

 It seems from Figure~\ref{Figure:ProfileNearEdge} that the weight of a point at a distance $d$ from the edge is of the form
 \[\left(1+h(d)\right)\frac{\vol(V)}{2\pi}\]
where  $h\colon (0,\infty)\to (-\infty,0)$ is a function such that
 \[h(d)\to -\infty \quad\text{as }d\to 0\qquad\text{and}\qquad h(d)\to 0\quad
 \text{exponentially as }d\to\infty.\]
 This behaviour is reminiscent of acutance in photographic images and the Gibbs phenomenon in signal processing.  However, I have not been able to figure out any heuristic giving a quantitative explanation of this.

\section{The bulk approximation and the penguin valuation}
\label{Section:Heuristic}
In this section we give a heuristic approximation to the contribution to the magnitude by the `bulk' of points in a `large' compact subset of $\R^n$.  Continuing in the heuristic vein and assuming some connection between magnitude and invariant valuations we are led to the definition of the penguin valuation $P$ at the end of this section.  In the next section we will see numerical evidence that the penguin valuation is often close to the magnitude.

In the previous section it was observed that in the `bulk' of a sufficiently large square all of the weights appeared to be constant.  Here we use that observation as the seed for a heuristic argument describing, for $X$ the closure of an open subset of $\R^n$, the contribution of the bulk of $X$ to the magnitude of $X$.  By $\bulk(X)$, the \emph{bulk} of $X$, we will mean all of the points `reasonably far from the boundary'.  Remember that the influence of the weight of a point decreases exponentially with distance, so you might take `reasonably far' to mean `at least ten units away': this is just a heuristic argument, after all.  The following is not a precise statement, but will be supported by numerical evidence in the next section.
\begin{approximation}
Let $X$ be the closure of a large open set in $\R^{n}$, let $\ddot X$ be the intersection of $X$ with a sufficiently fine lattice and let $\omega_{n}$ be the volume of the unit $n$-ball in $\R^{n}$.  Then the contribution to the magnitude of $\ddot X$ by the bulk of points in $\ddot X$ is roughly ${\vol (X)}/{n!\, \omega_{n}}$.
\end{approximation}
\begin{proof}[Heuristic justification]
The idea of the heuristic is that because the space is large and the effect of a point on other points decays exponentially with distance, the points in the bulk of $\ddot X$ only really know what is happening locally and so, as far as the weight equation is concerned, do not know that they are not part of an infinite lattice rather than a finite subset of one.  We can solve the weight equation for the infinite lattice and get a simple approximation to the resulting weight.  It is this approximation that we use to approximate the weight of points in the bulk of $\ddot X$.

So denote the lattice by $\LLL$.  We can solve the weight equations for $\LLL$, essentially using Speyer's method for homogeneous metric spaces (see~\cite[Theorem~1]{LeinsterWillerton:AsymptoticMagnitude}).  We will assume without loss of generality that the origin $0$ is in the lattice, so we can consider the weight equation for $0$:
 \[\sum_{x\in\LLL} e^{-d(x,0)}w_x =1.\]
We now assume, motivated by the observation in Section~\ref{Section:Squares} about the points in the bulk of the square, that all the points in $\LLL$ have the same weight $w$.  Writing, as usual, $|x|=d(0,x)$ we get
  \[w=\frac{1}{\sum_{x\in\LLL} e^{-|x|}}.\]
It is not immediate that the sum in the denominator converges, but we can approximate it by Riemann sums in the following way.

First pick a fundamental domain $V$ for $\LLL$;  for concreteness we can take this to be the Voronoi cell for the origin, which means all of the points in $\R^{n}$ which are at least as close to the origin as to any other lattice point.  The translates by $\LLL$ of the fundamental domain $V$ give a tessellation of $\R^{n}$, and each of the cells of the tessellation contains one point of the lattice.  The volume of the fundamental domain will be denoted by $\vol(V)$ and is sometimes called the \emph{covolume} of the lattice.

We can use this tessellation of the plane to calculate a Riemann sum approximation of an integral.  For the function $e^{-|x|}$, with the lattice sufficiently fine, we find
 \[\sum_{x\in\LLL} e^{-\left|x\right|}\vol(V)
 \simeq \int_{x\in \R^{n}}e^{-\left|x\right|} \,dx.\]
 The integral can be calculated by a change of variable to `polar coordinates':
 \[\int_{x\in \R^{n}}e^{-\left|x\right|} \,dx
   = \int_{\theta\in S^{n-1}} \int_{r=0}^{\infty}e^{-r}r^{n-1}\,dr\, d\theta
   = (n-1)!\,\vol(S^{n-1}),\]
 where the equality $\int_{r=0}^{\infty}e^{-r}r^{n-1}\,dr
   = (n-1)!$ can be shown using integration by parts and induction --- it is just the evaluation of a $\Gamma$ function.  As the the volume of the $(n-1)$-sphere and the $n$-ball are related by $\vol(S^{n-1})=n\vol(B^{n})=:n\omega_{n}$ we see
   \[\sum_{x\in\LLL} e^{-\left|x\right|}\vol(V)\simeq n! \,\omega_{n}\]
 and thus the weight of every point in the lattice $\LLL$ is given by
  \[w=\frac{\vol(V)}{\sum_{x\in\LLL} e^{-\left|x\right|}\vol(V)}
    \simeq \frac{\vol(V)}{n! \,\omega_{n}}.\]

Now return to the case of the finite metric space $\ddot X$ which is the intersection of the large space $X$ with a lattice $\LLL$.    As $X$ was taken to be the closure of an open set in $\R^n$, if we consider a point $x$ in the bulk, then as far $x$ is concerned it might as well be in the lattice $\LLL$ rather than $X\cap\LLL$, because the weight equation says that its interactions with points at a distance $r$ is $e^{-r}$, so it has exponentially small interaction with far-away points.  So assuming local homogeneity, we take the weight of $x$ to be the same as if it were considered as a point in $\LLL$ rather than $\ddot X$.  So the points in the bulk of $\ddot X$ all have roughly the same weight which is
\[w=\frac{\vol(V)}{n! \,\omega_{n}}.\]
The contribution of the points in the bulk to the magnitude is given by summing these weights, so we obtain roughly
 \[\frac{\vol(\bulk X)}{n! \,\omega_{n}} \simeq \frac{\vol( X)}{n! \,\omega_{n}},\]
as required.
\end{proof}

There are lots of approximations in the above argument, so it is not clear that it gives a sensible answer.  However, the answer is supported by the empirical data given in the next section; moreover, if we continue in a heuristic manner, as follows, we obtain an answer which is an even better fit with the empirical data.

The above heuristic indicates that for $X$ the closure of a large open subset of $\R^n$ then
  \[|X|\simeq \frac{\vol(X)}{n!\,\omega_n} +\text{contributions from near the boundary.}\]
Based on the rather meagre examples from~\cite{LeinsterWillerton:AsymptoticMagnitude}, we now assume that  for large $X$ the magnitude $|X|$ is close to some valuation $P(X)$.  Recall from~\cite{KlainRota:Book} and~\cite{LeinsterWillerton:AsymptoticMagnitude} that the intrinsic volumes $\{\mu_i\}$form a basis for the space of valuations.  As the $\mu_n$ is the $n$-dimensional volume of a subset of $\R^n$ we know that in dimension $n$ we must have $P(X)=\frac{\mu_n(X)}{n!\,\omega_n}+\text{lower order terms}$.  Now if $\iota\colon \R^i\hookrightarrow \R^n$ is an inclusion of Euclidean spaces and $Y\subset \R^i$ is some polyconvex subset then by the naturality of the intrinsic volumes with respect to inclusion of Euclidean spaces
 \begin{align*}P(\iota Y)&\simeq |\iota Y|=
   |Y|
   \simeq
   \frac{\mu_i(Y)}{i!\,\omega_i}+\text{lower order terms}\\
   & =
   \frac{\mu_i(\iota Y)}{i!\,\omega_i}+\text{lower order terms}.
\end{align*}
 So the coefficient of $\mu_i$ in $P$ is of the same form for all $i$.  We therefore define, for want of a better name, the \emph{penguin valuation} $P$ by
 \[P(X):= \sum_i \frac{\mu_i(X)}{i!\,\omega_i}\]
and speculate whether or not $|X|\simeq P(X)$ for large spaces $X$.

Again, this looks like rather wild speculation; however, when, as in the next section, for certain spaces $X$ we plot $P(tX)$ against the magnitude of approximations to $tX$ an impressively close fit is seen across a wide range of values of $t$.  This is especially true for the convex sets considered, namely the square, the disc and the cube, where $P(tX)$ is seen to be close to the approximation to $|tX|$ for all $t$, not just for large $t$.

From this, one is lead to ask if there is an interesting class of spaces $X$ for which the magnitude is asymptotically just the valuation $P$, i.e., $|tX|-P(tX)\to 0$ as $t\to\infty$, and is it true that for $X$ convex the magnitude is just the valuation, i.e., $|X|=P(X)$.

\section{Numerically comparing the magnitude with the penguin valuation}
\label{Section:Examples}
In this section we will plot the magnitude of various subsets of Euclidean spaces --- squares, discs, cubes, annuli and tori --- together with the penguin valuation of these subsets.   We can exactly calculate the penguin valuation for these and for the magnitude we numerically evaluate the magnitude of finite approximations to these subsets with around $25,000$ points.  We see a good correlation between the magnitude and the penguin valuation.  However, as pointed out in the introduction, this correlation seems unlikely to hold for every subset of Euclidean space.

 In each of the following cases we take an interesting subset $X$ of Euclidean space and approximate it by $\ddot X$ a finite subset of $X$ consisting of $N$ points.  We pick $N$ as large as possible subject to constraints of available computing power; these constraints are discussed in the appendix.  We then  use \texttt{maple} to calculate the magnitude $|t\ddot X|$ for various scale factors $t$; again, see the appendix for details of the \texttt{maple} code.
  We know from \cite{LeinsterWillerton:AsymptoticMagnitude} that for a finite metric space in which the points are far apart from one another, the magnitude will essentially be the number of points, this means that $|t\ddot X|\to N$ as $t\to \infty$.  However, for $t$ not so large, $|t\ddot X|$ should give a reasonable approximation to the mythical magnitude $|tX|$ of the space we are actually interested in.  Here, empirically, ``not so large''  seems to mean that the points are at most of order $0.1\,\text{units}$ from their nearest neighbour.

  Based on the ideas of the previous section we want to compare the magnitude $|tX|$ with the penguin valuation $P(tX)$ as we suspect  that they might be the same asymptotically, at least for when $X$ is the closure of an open subset.  So for $t$ large, but not too large, we suspect $|t\ddot X|\simeq |tX|$ and $|tX|\simeq P(tX)$, thus we would expect $|t\ddot X|\simeq P(tX)$ and it is this which is observed in the graphs below.  The graphs are plotted on log-log axes so that wide range of orders of magnitude can be observed.  It is striking that for the \emph{convex} examples we see that $|t\ddot X| \simeq P(tX)$ even for \emph{small} $t$, providing the motivation for the conjecture that $| X| = P(X)$ whenever $X$ is convex.

\subsection{Squares}
\begin{figure}[ht]
\begin{center}
\begin{tikzpicture}
\begin{loglogaxis}[tick label style={font=\footnotesize},
axis x line=bottom, axis y line = left,
xmin=0.08,ymin=0.8, ymax=60000, xmax=1000,
xtick={0.1,1,10,100,1000}, xticklabels={0.1,$1$,$10$,$100$,$1000$},
ytick={1,10,100,1000,10000,22500}, yticklabels={$1$,$10$,$100$,$1000$,$10000$,22500},
x axis line style={style = -},y axis line style={style = -},
xlabel={$t$}, x label style={at={(0.75,0)}},
legend style={at={(1,0.1)},anchor=south east,cells={anchor=west}}]
\addplot[only marks,mark=+] file {maple/SquareData.table};
\addplot [mark=none] expression[domain=0.1:800]  {1+x + x^2/6.2832};
\addplot[dashed,mark=none,color=blue] coordinates {(1000,22500) (0.08,22500)};
\legend{$|t\ddot Q|$,$P(tQ)$,$N$};
\end{loglogaxis}

\end{tikzpicture}

\end{center}
\caption{Squares: comparison of the penguin valuation $P(tQ)=t^2/2\pi + t+1$ with $|t\ddot Q|$, an approximation to the magnitude.}
\label{Figure:GraphSquare}
\end{figure}
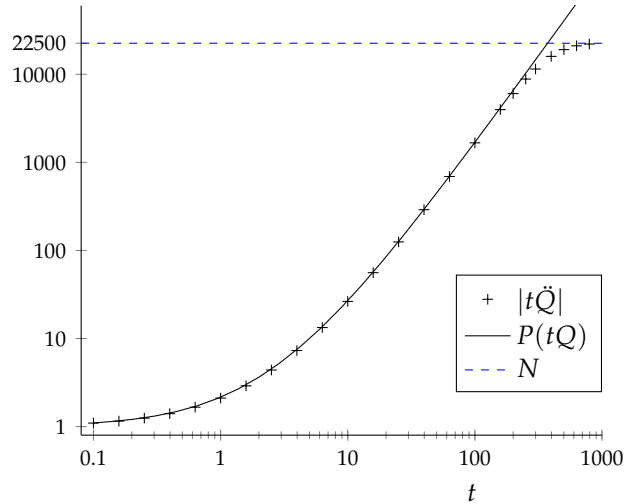

Here we consider  a $1\times 1$ square $Q$ which we approximate by a $150\times 150$ grid of points $\ddot Q$, this has total number of points, $N$, of $22,500$.  Of course the square $Q$ can be approximated by say a hexagonal grid or a large random set of points inside the square, but these appear to give similar results.  Weightings for squares with a square grid approximation were seen in  Figures~\ref{Fig:WeightsOfSquaresl2} and~\ref{Fig:SquaresVaryingWidths}.  We scale the finite square $\ddot Q$ by a factor of $t$ for various values of $t$ from $0.1$ to $1000$, numerically calculate the magnitude $|t\ddot Q|$ and plot the magnitudes against $t$ in Figure~\ref{Figure:GraphSquare}.  
We also plot $P(tQ)$ on the graph, which is given by
  \[P(tQ)=\frac{t^2}{2\pi}+\frac{2t}{2}+1,\]
as the linear intrinsic volume, $\mu_1$,  of a planar region is given by half the perimeter, which in this case is $2t$.
The curve gives a remarkable fit to the data points and fuels the suspicion that $|X|=P(X)$ when $X$ is a convex set.

\subsection{Discs}
\begin{figure}[th]
\begin{center}
\includegraphics{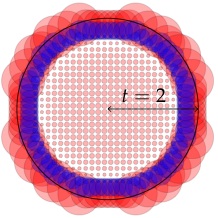}
%
\begin{tikzpicture}
\begin{loglogaxis}[tick label style={font=\footnotesize},
width = 0.6\textwidth,
axis x line=bottom, axis y line = left,
xmin=0.08,ymin=0.8, ymax=60000, xmax=1000,
xtick={0.1,1,10,100,1000}, xticklabels={0.1,$1$,$10$,$100$,$1000$},
ytick={1,10,100,1000,10000,25132}, yticklabels={$1$,$10$,$100$,$1000$,$10000$,25132},
x axis line style={style = -},y axis line style={style = -},
xlabel={$t$}, x label style={at={(0.75,0)}},
legend style={at={(1,0.1)},anchor=south east,cells={anchor=west}}]
\addplot[only marks,mark=+] file {maple/DiscData.table};
\addplot [mark=none] expression[domain=0.1:1000]  {1 + 1.570796327*x + 0.5*x^2};
\addplot[dashed,mark=none,color=blue] coordinates {(1000,25132) (0.08,25132)};
\legend{$|t\ddot D|$,$P(tD)$,$N$};

\end{loglogaxis}

\end{tikzpicture}
\end{center}
\caption{Discs: comparison of the penguin valuation $P(tD)=t^2/2 +\pi t/2+1$ with $|t\ddot D|$, an approximation to the magnitude.}
\label{Figure:GraphDisc}
\end{figure}
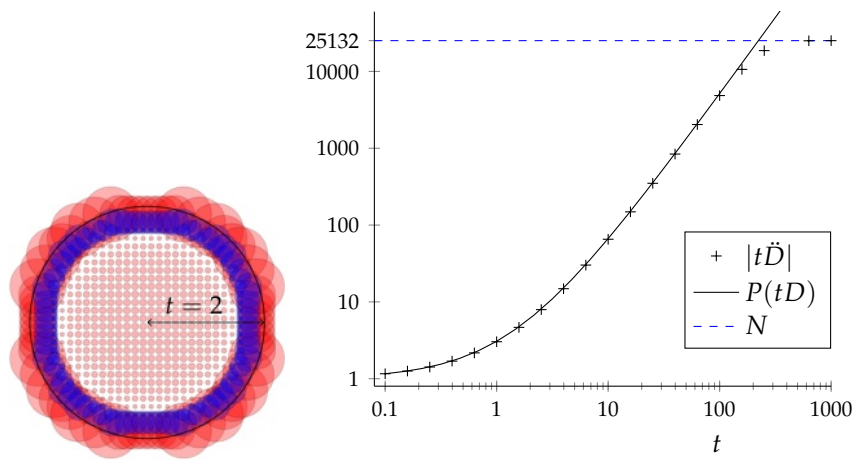

Similar to the last example, we consider $D$ to be a circle of radius $1$ and take $\ddot D$ to be a finite, square grid approximation to $D$; in this calculation we choose the mesh size of the grid to be such that we have approximately $25,000$ points in $D$.  We numerically calculate $|t\ddot D|$ the magnitude of the finite set $\ddot D$ scaled by a factor of $t$, taking various values of $t$ from $0.1$ to $1,000$, and in Figure~\ref{Figure:GraphDisc} we plot these values along with the graph of $P(tD)$, the suspected asymptotic magnitude of the disc $D$, where
 \[P(tD)=\frac{\pi t^2}{2\pi}+\frac{\pi t}{2}+1.\]
Again, we get a remarkably good fit with the data points up until the data points start flattening out.

\subsection{Cubes}
Again we consider a convex set, the cube, but this is one dimension higher than the above two examples of the square and the disc.  Let $C$ be the cube of side length $1$.  Now we approximate it by a $30\times 30 \times 30$ cubic lattice to get the total number of points, $N=27,000$, close to our computational maximum; this means that we get a less good approximation to the cube than we had for the square, as the points are five times further away from each other and, consequently, we see that the data points flatten out at smaller values of $t$ than for the square.  This is a problem of looking at higher dimensional sets.   However, there is still a very good fit with the graph of the penguin valuation
\[P(tC)=\frac{t^3}{8\pi} + \frac{3t^2}{2\pi}+\frac{3t}{2}+1,\]
for $t$ from $0.1$ to $40$.
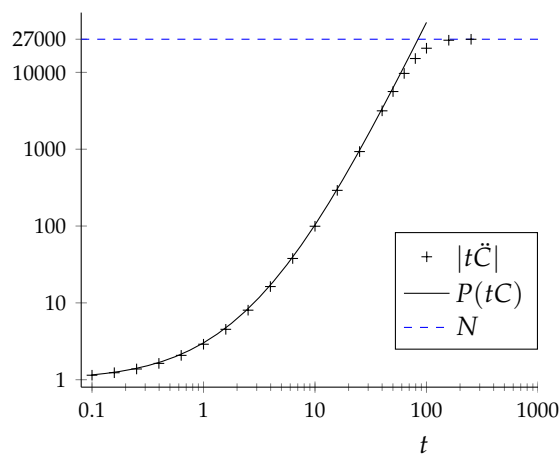
\begin{figure}
\begin{center}
\begin{tikzpicture}
\begin{loglogaxis}[tick label style={font=\footnotesize},
width = 0.6\textwidth,
axis x line=bottom, axis y line = left,
xmin=0.08,ymin=0.8, ymax=60000, xmax=1000,
xtick={0.1,1,10,100,1000}, xticklabels={0.1,$1$,$10$,$100$,$1000$},
ytick={1,10,100,1000,10000,27000}, yticklabels={$1$,$10$,$100$,$1000$,$10000$,27000},
x axis line style={style = -},y axis line style={style = -},
xlabel={$t$}, x label style={at={(0.75,0)}},
legend style={at={(1,0.1)},anchor=south east,cells={anchor=west}}]
\addplot[only marks,mark=+] file {maple/CubeData.table};
\addplot [mark=none] expression[domain=0.1:100]
    {1+1.5*x + 0.47746*x^2 +0.0397887*x^3};
\addplot[dashed,mark=none,color=blue] coordinates {(1000,27000) (0.08,27000)};
\legend{$|t\ddot C|$,$P(tC)$,$N$};
\end{loglogaxis}

\end{tikzpicture}

\end{center}
\caption{Cubes: comparison of the penguin valuation $P(tC)=t^3/8\pi+3t^2/2\pi +3 t/2+1$ with $|t\ddot C|$, an approximation to the magnitude.}
\label{Figure:GraphCube}
\end{figure}

\subsection{Annuli}
\begin{figure}[ht]
\begin{center}
\includegraphics{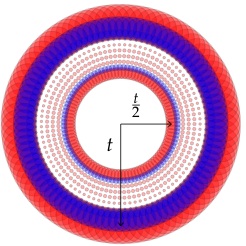}
%
%
\begin{tikzpicture}
\begin{loglogaxis}[tick label style={font=\footnotesize},
width = 0.6\textwidth,
axis x line=bottom, axis y line = left,
xmin=0.08,ymin=0.8, ymax=60000, xmax=1000,
xtick={0.1,1,10,100,1000}, xticklabels={0.1,$1$,$10$,$100$,$1000$},
ytick={1,10,100,1000,10000,23600}, yticklabels={$1$,$10$,$100$,$1000$,$10000$,23600},
x axis line style={style = -},y axis line style={style = -},
xlabel={$t$}, x label style={at={(0.75,0)}},
legend style={at={(1,0.1)},anchor=south east,cells={anchor=west}}]
\addplot[only marks,mark=+] file {maple/AnnulusData.table} node[right] {$N=30000$};
\addplot [mark=none] expression[domain=0.1:1000]  {2.356*x+0.375*x^2};
\addplot[dashed,mark=none,color=blue] coordinates {(1000,23600) (0.08,23600)};
\legend{$|t\ddot A|$,$P(tA)$,$N$};
\end{loglogaxis}

\end{tikzpicture}
\end{center}
\caption{Annuli: comparison of the penguin valuation $P(tA)=3t^2/8 +3\pi t/4$ with $|t\ddot A|$, an approximation to the magnitude.}
\label{Figure:GraphAnnulus}
\end{figure}
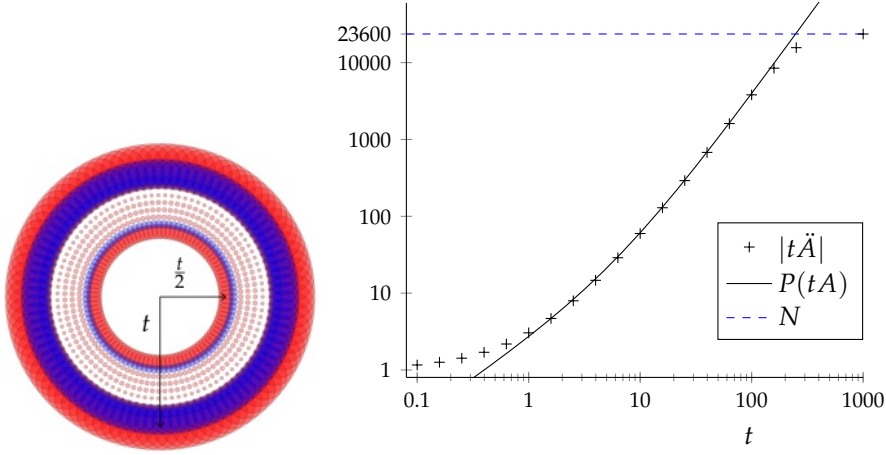

We now look at our first non-convex example.  We consider $A$, the annulus with outer radius $1$ and inner radius $1/2$.  This time we approximate not with a square grid but with $\ddot A$ a grid in polar coordinates.  In Figure~\ref{Figure:GraphAnnulus} we compare the calculated the calculated magnitudes with
 \[P(tA)=\frac{\frac34 \pi t^2}{2\pi}+\frac{\frac32\pi t}{2}.\]
In this case we don't see a good fit at small $t$, which we didn't necessarily expect anyway, but we certainly see a good fit for medium-sized $t$.

\subsection{Tori}
\begin{figure}[th]
\begin{center}
\begin{tikzpicture}
\begin{loglogaxis}[tick label style={font=\footnotesize},
axis x line=bottom, axis y line = left,
xmin=0.08,ymin=0.8, ymax=60000, xmax=1000,
xtick={0.1,1,10,100,1000}, xticklabels={0.1,$1$,$10$,$100$,$1000$},
ytick={1,10,100,1000,10000,25920}, yticklabels={$1$,$10$,$100$,$1000$,$10000$,25920},
x axis line style={style = -},y axis line style={style = -},
xlabel={$t$}, x label style={at={(0.75,0)}},
legend style={at={(1,0.1)},anchor=south east,cells={anchor=west}}]
\addplot[only marks,mark=+] file {maple/TorusData.table};
\addplot [mark=none] expression[domain=0.1:200]  {1.2566*x^2};
\addplot[dashed,mark=none,color=blue] coordinates {(1000,25920) (0.08,25920)};
\legend{$|t\ddot U|$,$P(tU)$,$N$};
\end{loglogaxis}

\end{tikzpicture}
\end{center}
\caption{Tori: comparison of the penguin valuation $P(tU)=t^2/10$ with $|t\ddot U|$, an approximation to the magnitude.}
\label{Figure:GraphTorus}
\end{figure}
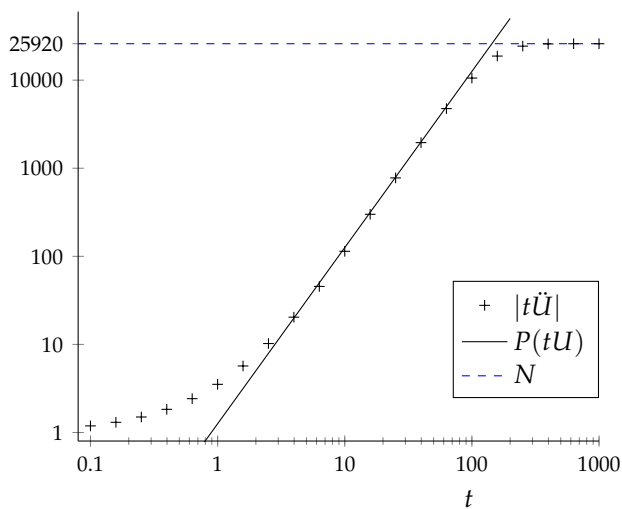

In this example we consider a subset of $\R^3$ which does not fit the hypotheses of the Bulk Approximation: the torus is a closed submanifold and is not the closure of an open subset of $\R^3$.  Nonetheless, the torus does seem to have the asymptotic behaviour predicted by the Bulk Approximation.  We consider $U$ to be the standardly embedded torus with major radius $1$ and minor radius $1/5$, in other words the locus of points
 \[\left( (1+\tfrac15\cos\phi) \cos\theta,(1+\tfrac15\cos\phi)\sin\theta,\tfrac15\sin\phi\right)
\quad\text{for }(\theta,\phi)\in[0,2\pi)\times[0,2\pi).\]
We get the finite approximation $\ddot U$ by taking $(\theta,\phi)$ to lie on a $m\times m$ grid of points in $[0,2\pi)\times[0,2\pi)$.  In Figure~\ref{Figure:GraphTorus} we compare the resulting magnitude data with the graph of the function
  \[P(tU)=\frac{\pi t^2/5}{2\pi}\]
and again see a good fit before the magnitudes of the finite space flatten out.

\section{Sierpinski gasket}
\label{Section:Sierpinski}
This example is of a different nature to the ones in Section~\ref{Section:Examples} so is given a section of its own.  In~\cite{LeinsterWillerton:AsymptoticMagnitude} we considered the fractal example of the ternary Cantor set.  It is not a polyconvex set and we can not make sensible interpretation of the intrinsic volumes $\{\mu_i\}_{i=0}^\infty$.  However we observed that the magnitude of  Cantor sets asymptotically satisfies an inclusion-exclusion principle in the following sense.  Suppose $P$ is some function on a set of subsets of Euclidean space which satisfies the inclusion-exclusion principle, so that $P(A)+P(B)-P(A\cap B)=P(A\cup B)$.  If $T_t$ denotes the Cantor set of length $t$, then we have the self-similarity homeomorphism for all $t>0$:
   \[T_{3t}\cong T_t\sqcup T_t.\]
 Thus if the Cantor sets are in the domain of $P$ and we set $p(t):=P(T_t)$ then
  \[p(3t)=2p(t).\]
The general solution to this functional equation is of the form $p(t)=f(t)t^{\log_32}$, where $f$ is a multiplicatively periodic function with $f(t)=f(3t)$ for all $t>0$.  In~\cite{LeinsterWillerton:AsymptoticMagnitude} we gave the explicit formula for such a function $p$ so that,  as $t\to\infty$, asymptotically, $|T_t|$ is $p(t)$; the corresponding function $f$ in that case is almost-constant, $f(t)\simeq 1.3$ for all $t>0$.

Now we can consider the Sierpinski gasket.  Unlike the case of the Cantor set we cannot calculate the magnitude exactly so we will have to numerically approximate it.  We can, however, carry out a similar analysis involving the inclusion-exclusion principle.  If $S_t$ denotes the Sierpinski gasket of side length $t$ then for $t>0$ the self-similarity can be represented by the homeomorphism
  \[S_{2t}\cong S_t  \sqcup S_t \sqcup S_t \smallsetminus \{3\text{ points}\},\]
so if $P$ a function satisfying the inclusion-exclusion principle, with the Sierpinski gaskets within its domain and with $P(\text{point})=1$, then defining $\widehat p(t):=P(S_t)$ we obtain the functional equation for all $t>0$:
  \[\widehat p(2t)=2\widehat p(t)-3.\]
The general solution to this is straightforwardly seen to be be $\widehat p(t)=\widehat f(t)t^{\log_23}+3/2$, where $\widehat f$ is a multiplicatively periodic function with $\widehat f(2t)=\widehat f(t)$ for all $t>0$.  Based on the behaviour of the Cantor sets, we might expect then that the magnitude $|S_t|$ of the Sierpinski gaskets asymptotically takes this form as $t\to \infty$.  We might also expect that the corresponding $\widehat f$ might is almost-constant.

\begin{figure}[ht]
\begin{center}
\begin{tikzpicture}[scale=0.8]
\input{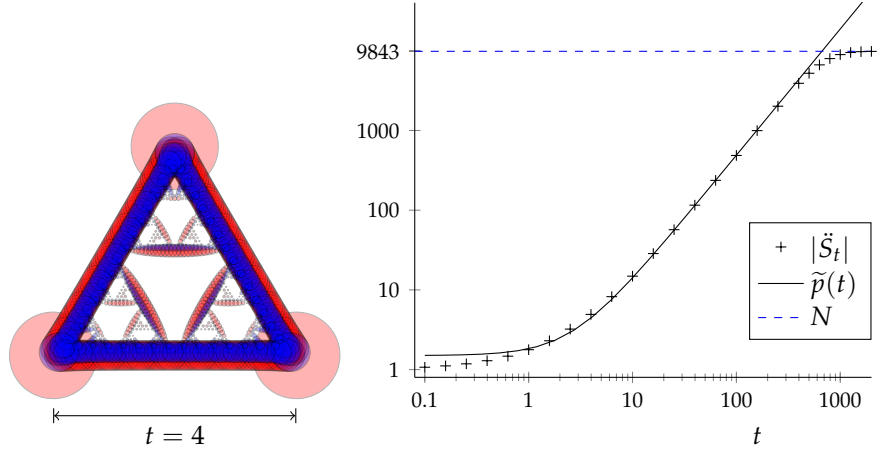}
\draw[|<->|] (0,-1.0)--(4,-1.0);
\draw (2,-1.0) node[below] {$t=4$};
\end{tikzpicture}
%
%
\begin{tikzpicture}
\begin{loglogaxis}[tick label style={font=\footnotesize},
width = 0.6\textwidth,
axis x line=bottom, axis y line = left,
xmin=0.08,ymin=0.8, ymax=40000, xmax=2000,
xtick={0.1,1,10,100,1000}, xticklabels={0.1,$1$,$10$,$100$,$1000$},
ytick={1,10,100,1000,9843,10000},
yticklabels={$1$,$10$,$100$,1000,$9843$,${}$},
x axis line style={style = -},y axis line style={style = -},
xlabel={$t$}, x label style={at={(0.75,0)}},
legend style={at={(1,0.1)},anchor=south east,cells={anchor=west}}]
\addplot[only marks,mark=+] file {maple/Sierpinski8Data.table};
\addplot [mark=none] expression[domain=0.1:2000]  {0.33*x^1.58496+1.5};
\addplot[dashed,mark=none,color=blue] coordinates {(2000,9843) (0.08,9843)};
\legend{$|\ddot S_t|$,$\widetilde p(t)$,$N$};
\end{loglogaxis}
\end{tikzpicture}
\end{center}
\caption{Sierpinski gasket: comparison of the calculated magnitudes $|\ddot S_t|$ with the function $\widetilde p(t):=\tfrac13 t^{\log_2 3}+\tfrac32$.}
\label{Figure:GraphSierpinski}
\end{figure}

For various values of $t$ the magnitude $|S_t|$ was approximated by numerically calculating the magnitude of $\ddot S_t$, a finite approximation to $S_t$ with $9843$ points.  The results are plotted in Figure~\ref{Figure:GraphSierpinski} where $\frac13 t^{\log_23}+\frac32$ is also plotted.  Here the function $\widehat f$ was taken to be the constant $1/3$ purely because that gave a good fit by eye to the calculated data.

The plot supports the guess that we have an asymptotic inclusion-exclusion principle for the magnitude, so there is an almost-constant function $\widehat f$ with $\widehat f(t)\simeq \frac13$ for all $t>0$, giving rise to the corresponding function $\widehat p$ as above with
  \[|S_t|-\widehat p(t)\to 0 \qquad\text{as }t\to\infty.\]
  
  One can observe that this is consistent with the more general guess that for $X$ a compact subset of $\R^n$, the asymptotic growth rate of the magnitude $|tX|$, i.e., 
    \[\lim_{t\to\infty}\frac{\log|tX|}{\log t},\]
is the Hausdorff dimension of $X$.
This is in agreement with the notion that for $X$ the closure of an open subset of $\R^n$ we have
 \[|tX|=\frac{t^n\vol(X)}{n!\,\omega_n}+\text{lower order terms}.\]
The next question to ask is what is the ``coefficient'' of the highest power of $t$ is in $|tX|$.  Whereas in some cases it appears to be proportional to the volume of $X$, in the fractal case it isn't even necessarily a constant but is an almost-constant function.  This is more food for thought.

\appendix
\section{Computer considerations}
The programming code involved in finding a weighting for a finite metric space is very simple as we are just solving $N$ linear equations in $N$ unknowns, where $N$ is the number of points.  Computer algebras have good algorithms built-in for solving such problems, so we can just rely on this.  This can be done as follows.  If $X$ is a metric space with points called $P_1,\dots, P_N$, then define $Z$ to be the matrix of exponentiated distances, define $w$ to be the vector of weights and define $\mathbf{1}$ to be the vector of ones, so 
\[
  Z_{ij}:=\exp\bigl(-d(P_i,P_j)\bigr); \quad 
  w:=(w_{P_1},\dots,w_{P_N})^\text{T}\in \R^N; \quad
  \mathbf{1}:=(1,\dots,1)^\text{T}\in \R^N.
\]
The weight equations are then expressed as the single vector equation
\[Zw=\mathbf{1}.\]
We need to solve this for $w$ and add up the entries of $w$ to obtain the magnitude of $X$.  Here is a piece of pseudo \texttt{maple} code, so you can see how straightforward it is.

\begin{code}

 # Load in a library so we can use Matrix, Vector and LinearSolve
with(LinearAlgebra):

 # Define N to be the number of points

 # Define P to be a list with N entries so that
 # P[i]:=[i-th x-coordinate, i-th y-coordinate]:
 # or possibly with a z-coordinate as well

 # Define the distance d(a,b) between
 # points a=[a[1],a[2]] and b=[b[1],b[2]],
 # adding a third term if working in 3-space
d:=(a,b)->sqrt((a[1]-b[1])^2+(a[2]-b[2])^2):

 # Define the matrix of exponentiated distances
Z:=Matrix(1..N,1..N,
           (i,j) -> exp(-d(P[i],P[j])),
           datatype=float, shape=symmetric):

 # Solve the weight equations and put the solution in the vector w
 # so that w[i] is the weight of the i-th point
w:=LinearSolve(Z,Vector(1..N,1));

 # Add all of the weights together to give the magnitude
Magnitude:=add(w[i],i=1..N):

\end{code}

To get a good approximation we want to take the number of points, $N$, to be as large as possible.  The principal constraints on how large we can take $N$ to be are those of computing time and available memory.  I was using the \texttt{iceberg} cluster at the University of Sheffield, each processor there has 12GB RAM and that turned out to be the main constraint.  Taking $N$  to be around $25,000$ resulted in between 3GB and 12GB RAM being used, and between 5 and 12 hours of CPU time for the calculation of a single magnitude.  Each data point plotted in Section~\ref{Section:Examples} required that level of resources.  Interestingly the metric spaces which needed the most time and memory were those near the `shoulder' in each of the graphs of Section~\ref{Section:Examples} where the data-points begin to flatten out at the top.

I should add that I have given no thought to how accurate the answer \texttt{maple} gives is, in terms of what sort of rounding errors would be present after row reducing a $25,000\times 25,000$ matrix, but, nonetheless, the numerical answers obtained looked sufficiently convincing.
\bigskip

\end{document}

%% file: SWi.tex
\usepackage{amssymb}
\usepackage{mathrsfs}

\newcommand{\R}{\mathbb R}

\newcommand{\LLL}{\mathscr L}

\DeclareMathOperator{\vol}{vol}
\DeclareMathOperator{\bulk}{bulk}

\renewcommand{\paragraph}[1]{\bigskip\noindent \emph{#1.}}

\newtheorem*{approximation}{Bulk Approximation}

\theoremstyle{definition}

\renewcommand{\phi}{\varphi}
\renewcommand{\epsilon}{\varepsilon}